\documentclass[11pt,leqno,fleqn]{article}
\usepackage{epsf,amsmath,amssymb,amscd,hyperref,srcltx}
\usepackage{epsfcenter}
\usepackage{amsfonts}
\usepackage{counters}
\usepackage{psfrag}
\newcommand{\dps}{\displaystyle}

\newcommand{\ii}{\text{\rm i}}

\newcommand{\T}{^{\sf T}}
\newcounter{figuren}

\newcommand{\dy}[2]{%
\refstepcounter{equation}%
\LABEL{#1}%
\begin{list}{}{
\topsep 5mm
\leftmargin 18mm
\rightmargin 0cm
\itemsep 0mm
\listparindent 0mm
\parsep 0mm
\itemsep 0mm
\labelsep 0mm
\labelwidth 18mm
}%
\item[\rm (\theequation)\hfill]
#2
\end{list}%
}
\newcommand{\dyz}[1]{%
\refstepcounter{equation}%
\begin{list}{}{
\topsep 5mm
\leftmargin 18mm
\rightmargin 0cm
\itemsep 0mm
\listparindent 0mm
\parsep 0mm
\itemsep 0mm
\labelsep 0mm
\labelwidth 18mm
}%
\item[\rm (\theequation)\hfill]
#1
\end{list}%
}
\newcommand{\dyyz}[1]{\dyz{\raggedright$\dps#1$}}
\newcommand{\dyy}[2]{\dy{#1}{\raggedright$\dps#2$}}
\newcommand{\de}[2]{\dy{#1}{\raggedright$\displaystyle #2 $}}
\newcommand{\dez}[1]{\dyz{\raggedright$\displaystyle #1 $}}
\newcommand{\leeg}[1]{}

\newcommand{\di}[2]{%
\refstepcounter{equation}%
\LABEL{#1}%
\begin{list}{}{
\topsep 5mm
\leftmargin 10mm
\rightmargin 0cm
\itemsep 0mm
\listparindent 0mm
\parsep 0mm
\labelsep 1mm
\labelwidth 10mm
}%
\item[\rm (\theequation)\hfill]
\begin{list}{}{
\topsep 0mm
\leftmargin 8mm
\rightmargin 0mm
\itemsep 0mm
\listparindent 0mm
\parsep 0mm
\labelsep 1.5mm
\labelwidth 6.5mm
}
#2
\end{list}%
\end{list}%
}
\newcounter{stelling}
\newcommand{\thm}[2]{\setcounter{gevolg}{0}\setcounter{claim}{0}\refstepcounter{stelling}\vspace{4mm}\noindent{\bf Theorem \thestelling.}\label{#1}{\it #2}}

\newalphacounterwithin{gevolg}{stelling}
\newcommand{\cor}[2]{\refstepcounter{gevolg}\setcounter{claim}{0}\vspace{4mm}\noindent{\bf Corollary \thegevolg.}\label{#1}{\it #2}}

\newcounter{hulpstelling}

\newcounter{bewering}

\newcounter{claim}

\newalphacounterwithin{subclaim}{claim}

\newcounter{opmerking}

\newcounter{hoofdstuk}

\newcounter{sectie}
\newcounter{subsectie}
\newcounterwithin{ex}{sectie}

\newcommand{\sectz}[1]{\refstepcounter{sectie}\setcounter{subsectie}{0}\setcounter{ex}{0}
\section*{\boldmath \thesectie. #1}%
}

\newcounter{lit}

\newcommand{\pf}{\vspace{3mm}\noindent{\bf Proof.}\ }

\newcommand{\bx}{\hspace*{\fill} \hbox{\hskip 1pt \vrule width 4pt height 8pt depth 1.5pt \hskip 1pt}

\addvspace{4mm}}

\newcommand{\bxx}{\hspace*{\fill} \hbox{\hskip 1pt \vrule width 4pt height 8pt depth 1.5pt \hskip 1pt}}

\newcommand{\rf}[1]{{\rm (\ref{#1})}}

\newcommand{\CC}{{\cal C}}

\newcommand{\GG}{{\cal G}}

\newcommand{\LL}{{\cal L}}
\newcommand{\MM}{{\cal M}}

\newcommand{\OO}{{\cal O}}

\renewcommand{\SS}{{\cal S}}
\newcommand{\TT}{{\cal T}}

\newcommand{\VV}{{\cal V}}
\newcommand{\WW}{{\cal W}}

\newcommand{\kfrac}[2]{\mbox{$\frac{#1}{#2}$}}
\newcommand{\kint}[2]{\mbox{$\int$}}

\newcommand{\NIET}[1]{}

\newcommand{\LABEL}[1]{\label{#1}}

\newcommand{\sgn}{\text{\rm sgn}}

\newcommand{\oC}{{\mathbb{C}}}

\newcommand{\oN}{{\mathbb{N}}}

\newcommand{\oR}{{\mathbb{R}}}

\renewcommand{\Im}{{\text{\rm Im~}}}
\newcommand{\Ker}{{\text{\rm Ker~}}}

\newcommand{\MP}[1]{(\oC^{n\times n})^{\otimes #1}}

\begin{document}

\begin{center}
{\LARGE\bf On virtual link invariants

}

\end{center}
\vspace{1mm}
\begin{center}
{\large
\hspace{10mm}
Alexander Schrijver\footnote{ CWI and University of Amsterdam.
Mailing address: CWI, Science Park 123, 1098 XG Amsterdam,
The Netherlands.
Email: lex@cwi.nl.}}

\end{center}

\noindent
{\small
{\bf Abstract.}
We characterize the virtual link invariants that are
partition functions of vertex models (as considered by
de la Harpe and Jones), both in the real and in the complex case.
We show that for any fixed number of states, these invariants
form an affine variety.
Basic techniques are the first and second
fundamental theorem of invariant theory
for the orthogonal group (in the sense of Weyl)
and some related methods from algebraic geometry.

}

\sectz{Introduction and survey of results}

This paper is inspired by some recent results in the range of
combinatorial parameters and invariant theory,
in particular results of
Szegedy [12],
Freedman, Lov\'asz, and Schrijver [2],
and
Schrijver [10,\linebreak[0]11].
The concepts of virtual link and virtual link diagram were introduced
by Kauffman [6]; see Kauffman [7] for more background,
which we partially repeat below.

A {\em virtual link diagram} is an undirected graph $G$ such that
for each vertex $v$:
\di{18ok07a}{
$v$ has degree 4 and the edges incident with $v$ are cyclically
ordered, where one pair of opposite edges is labeled as
`over-going'.
}
$G$ may have loops and multiple edges.
Moreover, `vertexless' loops are allowed, that is, loops without a vertex.
We denote the vertexless loop by $O$.
Let $\GG$ denote the collection of virtual link diagrams,
two of them being the same if they are isomorphic.

A virtual link diagram can be seen as the projection
of a link in $M\times\oR$ where $M$ is some oriented surface.
Since this connection however is not stable under all Reidemeister
moves (e.g., one may need to create a handle to allow a
type II Reidemeister move), we will view virtual link diagrams
just abstractly as given above.

In this paper, $\oN=\{0,1,2,\ldots\}$ and
\dez{
[n]:=\{1,\ldots,n\}
}
for any $n\in\oN$.
We denote the sets of vertices and edges of a virtual
link diagram $G$ by $VG$ and $EG$, respectively.
$K_0$ denotes the virtual link diagram with no vertices and edges.

Let $V$ be an $n$-dimensional complex linear space with a symmetric
nondegenerate bilinear form $\langle.,.\rangle$.
We can identify $V$ with $\oC^n$,
with the standard bilinear form $\langle x,y\rangle=x\T y$.
Having the bilinear form, we can identify $V^*$ with $V$.

Let $S_2$ act on $V^{\otimes 4}$ so that the nonidentity element of
$S_2$ brings $x_1\otimes x_2\otimes x_3\otimes x_4$
to $x_3\otimes x_4\otimes x_1\otimes x_2$.
Define
\dez{
\MM_n:=(V^{\otimes 4})^{S_2},
}
that is, the linear space of $S_2$-invariant elements of $V^{\otimes 4}$.
Note that $\MM_n$ can be identified with the collection of symmetric
matrices in $(\oC^{n\times n})^{\otimes 2}$.

Following de la Harpe and Jones [4], any element $R$ of $\MM_n$ can be called
a {\em vertex model}.
For any $R\in\MM_n$,
let $f_R$ be the {\em partition function} of $R$;
that is, the virtual link diagram invariant defined by
\de{23ap07c}{
f_R(G)=\sum_{\phi:EG\to[n]}
\prod_{v\in VG}R_{\phi(\delta(v))}.
}
Here we put
\dez{
\phi(\delta(v)):=(\phi(e_1),\phi(e_2),\phi(e_3),\phi(e_4)),
}
where $e_1,e_2,e_3,e_4$ are the edges incident with $v$, in clockwise order,
and where $e_1,e_3$ is the over-going pair.
Since $R$ is $S_2$-invariant, $R_{\phi(\delta(v))}$ is well-defined.

We call $f:\GG\to\oC$ an {\em $n$-state partition function}
if $f=f_R$ for some $R\in\MM_n$.
In Corollary \ref{24no08b} we characterize $n$-state partition functions.

\sectz{Reidemeister moves}

Reidemeister moves yield an isotopy of virtual link diagrams.
A {\em virtual link} is an isotopy class of virtual
link diagrams.
A {\em virtual link diagram invariant} is a function defined
on virtual link diagrams that is invariant under
isomorphisms.
A {\em virtual link invariant} is a virtual 
link diagram invariant that is invariant under Reidemeister moves.
(So in fact it is a function on virtual links, but the
definition as given turns out to be more convenient.)

We recall the well-known sufficient conditions for $f_R$
to be a virtual link invariant (i.e., to be invariant under
Reidemeister moves).
To this end, let $C:\MP{2}\to\oC^{n\times n}$,
$D:\MP{2}\to\MP{2}$, and
$E_{1,2},E_{1,3},E_{2,3}:\MP{2}\to\MP{3}$ be
the unique linear functions satisfying
\dyz{
$C(M\otimes N)=MN$,
$D(M\otimes N)=M\otimes N\T$,
$E_{1,2}(M\otimes N)=M\otimes N\otimes I_n$,
$E_{1,3}(M\otimes N)=M\otimes I_n\otimes N$,
$E_{2,3}(M\otimes N)=I_n\otimes M\otimes N$
}
for all $M,N\in\oC^{n\times n}$,
where $I_n$ denotes the identity matrix in $\oC^{n\times n}$.
Then a {\em sufficient} condition for $f_R$ to be a
virtual link invariant is
(cf.\ Turaev [13], Kauffman [5]):
\dy{22ok07a}{
$C(R)=I_n$,
$RD(R)=I_n\otimes I_n$,
$E_{1,2}(R)E_{1,3}(R)E_{2,3}(R)=
E_{2,3}(R)E_{1,3}(R)E_{1,2}(R)$.
}
The last equation in \rf{22ok07a} is the {\em Yang-Baxter equation}.

In the real case, condition \rf{22ok07a} is necessary and sufficient
for $f_R$ being invariant under Reidemeister moves.
But in the complex case,
condition \rf{22ok07a} is only a sufficient, and not a necessary
condition for $f_R$ being invariant under Reidemeister moves.
For instance, if $R=A\otimes A$ with
$A:=$\mbox{\tiny $\left(\begin{array}{cc}2&\ii\\ \ii&0\end{array}\right)$}
(where $\ii$ is the imaginary unit),
then $f_R(G)$ is equal to $2^k$ where $k$ is the number of knots in the
virtual link $G$.
(A {\em knot} in a virtual link diagram $(V,E)$ is a component of the
graph on $E$ where two edges are adjacent if and only if there is a vertex
where they are oppositie.)
So $f_R$ is invariant under Reidemeister moves.
Yet, $C(R)\neq I_2$.

We will however see that there exists for each $R$ with $f_R$ invariant
under Reidemeister moves an $R'$ with $f_{R'}=f_R$ and $R'$ satisfying
\rf{22ok07a}.
This holds in fact more generally for any linear combination of
`tangles' that leaves $f_R$ invariant.

\sectz{Some notation and terminology}

Let $\oC\GG$ denote the space of formal linear combinations of
elements of $\GG$.
The elements of $\oC\GG$ are called {\em quantum virtual link diagrams}.
Any function on $\GG$ with values in a $\oC$-linear space can be
extended linearly to $\oC\GG$.
Let $GH$ denote the disjoint union of $G$ and $H$.
Then, taking $GH$ as products, makes $\oC\GG$ to an algebra.

As usual, let $\OO(X)$ denote the set of all regular functions
on a linear space $X$.
Define $p_n:\GG\to \OO(\MM_n)$ by
\dez{
p_n(G)(X):=f_X(G)
}
for $G\in\GG$ and $X\in\MM_n$.
Note that $p_n(K_0)=1$ and
\dez{
p_n(G)p_n(H)=p_n(GH).
}
Extending $p_n$ to $\oC\GG$,
makes $p_n$ to an algebra homomorphisms $\oC\GG\to\OO(\MM_n)$.

A central step in our proof is describing
the image $\Im p_n$ and the kernel $\Ker p_n$ of $p_n$.
This is based on the First and Second Fundamental Theorem of
Invariant Theory (in the sense of Weyl [14]) for the
orthogonal group (the `FFT' and `SFT') --- cf.\ Goodman and Wallach [3].

\sectz{Tangles}

Let $k\in\oN$.
A {\em $k$-tangle} is an undirected graph where each vertex
$v$ either satisfies \rf{18ok07a} or has degree 1, and where
there are precisely $k$ vertices of degree 1, labeled $1,\ldots,k$.
(It follows that $k$ is even.)
Let $\TT_k$ and $\TT$ denote the collections of $k$-tangles and of
tangles, respectively.
A {\em tangle} is a $k$-tangle for any $k$.
The $0$-tangles are precisely the virtual link diagrams, and
so $\TT_0=\GG$.

A {\em quantum $k$-tangle} is a formal $\oC$-linear combination of
$k$-tangles, and a {\em quantum tangle} a formal linear combination
of tangles.
Let $\oC\TT_k$ and $\oC\TT$ denote the spaces of quantum
$k$-tangles and quantum tangles, respectively.
(So $\oC\TT_0=\oC\GG$.)
Any function on the set of tangles with values
in a $\oC$-linear space can be linearly extended to $\oC\TT$.

Let $u$ and $v$ be vertices of degree 1 in a tangle.
We say that we {\em glue} $u$ and $v$ if we identify $u$ and $v$
and ignore them as vertices, joining the incident edges to one edge.
If $T$ and $T'$ are $k$-tangles, then $T\cdot T'$ is defined
to be the $0$-tangle obtained from the disjoint union of $T$ and $T'$ by
glueing the vertices labeled $i$ in $T$ and $T'$, for $i=1,\ldots,k$.
By defining $T\cdot T'=0$ for any $k$-tangle $T$ and
$l$-tangle $T'$ with $k\neq l$, we obtain a bilinear form
$\oC\TT\times\oC\TT\to\oC\GG$.

\sectz{Applying invariant theory}

As usually, for any group $G$ acting on a set $S$, let
$S^G:=\{x\in S\mid x^U=x$ for each $U\in G\}$.

Let $O_n$ denote the group of complex orthogonal matrices.
$O_n$ acts on $V^{\otimes k}$ by $x^U= U^{\otimes k}x$ for
$x\in V^{\otimes k}$.
In particular, $O_n$ acts on $\MM_n$.
This induces an action on $\OO(\MM_n)$.
Note that
\de{18de07ax}{
f_{R^U}(G)=f_R(G)
}
for each virtual link diagram $G$ and each $U\in O_n$.
Hence $p_n(G)$ is an $O_n$-invariant polynomial for each $G$.
As we will see, these span all $O_n$-invariant polynomials
in $\OO(\MM_n)$.

For each $\pi\in S_n$, let $T_{\pi}$ be the $2n$-tangle
with vertex set $[2n]$, labeled by $\lambda:[2n]\to[2n]$ with
$\lambda(i)=i$ for $i\in[2n]$,
and edges $\{i,n+\pi(i)\}$ for each $i\in[n]$.
Define the quantum $2n$-tangle ${\det}_n$ by
\dez{
{\det}_n:=\sum_{\pi\in S_n}\sgn(\pi)T_{\pi}.
}

Define
\dez{
J_n:={\det}_n\cdot\oC\TT_{2n}.
}

\thm{7ja09a}{
$\dps \Im p_n=\OO(\MM_n)^{O_n}$ and $\dps\Ker p_n=J_{n+1}$.
}

\pf
It is easy to see that each $p_n(G)$ is $O_n$-invariant.
This gives $\subseteq$ in the first equality.
One also checks directly $\supseteq$ in the second equality.

To see the reverse inclusions,
choose $k\in\oN$, and let $\GG_k$ be the collection
of virtual link diagram with $k$ vertices.
Let $\OO_k(\MM_n)$ be the set of all homogeneous polynomials 
in $\OO(\MM_n)$ of degree $k$.

Let $m:=4k$.
Let $S\oC^{m\times m}$ be the set of symmetric matrices in
$\oC^{m\times m}$.
Let $\LL$ be the subspace of $\OO(S\oC^{m\times m})$ spanned
by the monomials $\prod_{ij\in M}x_{ij}$ where $M$ is any perfect
matching on $[m]$.
We give each $j\in[m]$ label $j$.
This makes any perfect matching $M$ on $[m]$ to an $m$-tangle.

We define linear functions $\tau$, $\mu$, and $\sigma$
so as to make a commutative diagram:
\de{27no10a}{
\begin{CD}
\oC\GG_k @>p>> \OO_k(\MM_n)\\
@AA\mu A       @AA\sigma A\\ 
\LL @>\tau>>  (\oC^n)^{\otimes m}
\end{CD}
}
First, define $\mu$ and $\tau$ by
\dyz{
$\dps\mu(\prod_{ij\in M}x_{ij}):=T\cdot M$
~~~~and~~~~
$\dps\tau(\prod_{ij\in M}x_{ij}):=\sum_{\phi:[m]\to[n]\atop
\forall e\in M: |\phi(e)|=1}e_{\phi}$,
}
for any perfect matching $M$ on $[m]$.
Here $T$ is the $m$-tangle with labeled vertices $1,\ldots,m$
and unlabeled vertices $v_1,\ldots,v_k$, with $v_i$ connected
to vertices $4i-3,3i-2,4i-1,4i$, in this order.
Moreover,
\dez{
e_{\phi}:=\bigotimes_{h\in [m]}e_{\phi(h)}.
}

Next, define $\sigma$ by
\dyz{
$\dps\sigma(e_{\phi}):=\prod_{i\in[k]}
\rho_{S_2}(e_{\phi_i})
$
}
for $\phi:[m]\to[n]$.
Here $\rho_{S_2}$ is the Reynolds operator for $S_2$, and
$\phi_i$ is the function $[4]\to [n]$ with $\phi_i(j):=\phi(4i-4+j)$
for $j\in[4]$.

Now diagram \rf{27no10a} commutes, that is,
\de{26no10a}{
p\circ\mu=\sigma\circ\tau.
}
To prove it, choose a perfect matching $M$ on $[m]$.
Then $p(\mu(\prod_{ij\in M}x_{ij}))=p(T\cdot M)$.
Moreover,
\dyyz{
\sigma(\tau(\prod_{ij\in M}x_{ij}))=
\sum_{\phi:[m]\to[n]\atop
\forall e\in M: |\phi(e)|=1}\sigma(e_{\phi})
=
\sum_{\phi:[m]\to[n]\atop
\forall e\in M: |\phi(e)|=1}
\prod_{i\in[k]}
\rho_{S_2}(e_{\phi_i})
=
p(T\cdot M).
}
This proves \rf{26no10a}.

Now by the FFT, $\Im\tau=((\oC^n)^{\otimes m})^{O_n}$.
Moreover, $\Im\mu=\oC\GG_k$ and $\Im\sigma=\OO_k(\MM_n)$.
Hence
\dyyz{
p(\oC\GG_k)=\Im(p\circ\mu)=\Im(\sigma\circ\tau)=
\sigma(\Im\tau)=\sigma((\oC^n)^{\otimes m})^{O_n}
=(\OO_k(\MM_n))^{O_n}.
}
So $\Im p=\OO(\MM_n)^{O_n}$.

By the SFT, $\Ker\tau$ is equal to the ideal $I$ in
$\OO(S\oC^{m\times m})$ generated by the $(n+1)\times(n+1)$
minors of $S\oC^{m\times m}$.
Choose $\gamma\in\oC\GG_k$ with $p(\gamma)=0$.
Then $\gamma=\mu(q)$ for some $q\in\LL^{\SS}$,
where $\SS$ is the group of permutations of $[m]$ generated by $S_k$
and by $S_2$.
Hence $\tau(q)\in((\oC^n)^{\otimes m})^{\SS}$.
As $\sigma(\tau(q))=p(\mu(q))=p(\gamma)=0$, this implies
$\tau(q)=0$.
Hence $q\in I$.
Therefore, $\gamma=\mu(q)\in\mu(I)\subseteq J_{n+1}$.
\bx

\cor{24no08b}{
Let $f$ be a virtual link diagram invariant and let $n\in\oN$.
Then $f$ is an $n$-state partition function if and only if $f$ is multiplicative and
$f({\det}_{n+1}\cdot T)=0$ for each $2(n+1)$-tangle $T$.
}

\pf
Necessity follows from the fact that if $f=f_R$, then
$f(G)=p_n(G)(R)$ for each $G$.
So $f$ is multiplicative as $p_n$ is an algebra homomorphism.
Moreover, for any $2(n+1)$-tangle $T$ one has
$f_R({\det}_n\cdot T)=p_n({\det}_n\cdot T)(R)=0$,
as ${\det}_n\cdot T$ belongs to $\Ker p_n$.

To see sufficiency, as $f(\Ker p_n)=0$, there exists an algebra
homomorphism $\hat f:\OO(\MM_n)^{O_n}\to\oC$ such that
$\hat f\circ p_n=f$.
Now $1\not\in\Ker\hat f$.
Hence, since $O_n$ is reductive, $\Ker\hat f$ has a common zero $R$.
Then for each virtual link diagram $G$,
$f_R(G)-f(G)=(p_n(G)-f(G))(R)=0$, as $p_n(G)-f(G)$
belongs to $\Ker\hat f$.
\bx

\sectz{Extension to tangles}

Let $R\in\MM_n$ and let $T$ be a $k$-tangle.
Define the tensor $f_R(T)$ in $V^{\otimes k}$ by
\dez{
f_R(T):=
\sum_{\kappa:ET\to [n]}
\left(
\prod_{v\in V'T}R_{\kappa(\delta(v))}
\right)
e_{\kappa(\varepsilon_1)}\otimes\cdots\otimes e_{\kappa(\varepsilon_k)}.
}
Here $\varepsilon_i$ denotes the edge incident with the vertex labeled
$i$, for $i=1,\ldots,k$.
Moreover,
\dyz{
$V'T:=$ set of unlabeled vertices of $T$.
}

Define
\dez{
T(V):=\bigoplus_{k=0}^{\infty}V^{\otimes k}.
}

The set $\OO(\MM_n)\otimes V^{\otimes k}$ is naturally
equivalent to the set of all morphisms $\MM_n\to V^{\otimes k}$.
(The element $p\otimes v$ with $p\in\OO(\MM_n)$ and $v\in V^{\otimes k}$
corresponds to the function $x\mapsto p(x)v$.)
Then $O_n$-invariant elements in $\OO(\MM_n)\otimes T(V)$
correspond to the $O_n$-equi\-va\-ri\-ant morphisms
$\MM_n\to T(V)$.

The bilinear form on $V$ induces a bilinear form
$\cdot:T(V)^2\to\oC$, which induces a
bilinear form
$\cdot:(\OO(\MM_n)\otimes T(V))^2\to \OO(\MM_n)$.
It brings pairs of $O_n$-equivariant morphisms to
$O_n$-invariant polynomials in $\OO(\MM_n)$.

For any $k$-tangle $T$, define the element $p_n(T)$ of
$\OO(\MM_n)\otimes V^{\otimes k}$ by
\dez{
p_n(T)(X):=f_X(T)
}
for $X\in\MM_n$.
Note that this extends the definitions of $p_n(G)$ and $f_R(G)$
given above.
Then for quantum tangles $\tau,\tau'$:
\de{24ok07b}{
p_n(\tau\cdot\tau')=p_n(\tau)\cdot p_n(\tau').
}

\thm{2no07b}{
$p_n(\oC\TT)=(\OO(\MM_n)\otimes T(V))^{O_n}$.
}

\pf
This can be proved similarly to Theorem \ref{7ja09a}, now taking
$m=4k+l$.
\bx

\sectz{Nondegeneracy}

Define
\dyz{
$\WW_n:=\{R\in\MM_n\mid$
the bilinear form $\cdot$ on $f_R(\oC\TT)$ is nondegenerate$\}$.
}
The condition is equivalent to saying that
for each quantum tangle $\tau$: if $f_R(\tau\cdot T')=0$
for each tangle $T'$, then $f_R(\tau)=0$.

\thm{23no07b}{
For each $R\in\MM_n$ there exists $R'\in\WW_n$ with $f_{R'}(G)=f_R(G)$
for each virtual link diagram $G$.
}

\pf
Let $\CC$ be the collection of all quantum tangles $\tau$ such that
$f_R(\tau\cdot T')=0$ for each tangle $T'$.
We first show that the $p_n(\tau)$ (over all $\tau\in\CC$)
have a common zero.

Suppose such a common zero does not exist.
Then by the Nullstellensatz there exists a finite subset $\CC_0$ of
$\CC$ and for each $\tau\in\CC_0$ a $q_{\tau}\in\OO(\MM_n)\otimes T(V)$
with
$\sum_{\tau\in\CC_0}p_n(\tau)\cdot q_{\tau}=1$.
Applying the Reynolds operator we can assume that
each $q_{\tau}$ belongs to $(\OO(\MM_n)\otimes T(V))^{O_n}$.
So by Theorem \ref{2no07b},
$q_{\tau}=p_n(\tau')$ for some tangle $\tau'$ (depending on $\tau$).
Then
\dyy{9de10a}{
1=
\sum_{\tau\in\CC_0} p_n(\tau)\cdot p_n(\tau')
=
\sum_{\tau\in\CC_0} p_n(\tau\cdot \tau').
}
However, $p_n(\tau\cdot\tau')(R)=0$ for each $\tau\in\CC$,
since $f_R(\tau\cdot\tau')=0$.
So \rf{9de10a} does not hold.

So the $p_n(\tau)$ have a common zero $R'$.
Then $f_{R'}(G)=f_R(G)$ for any virtual link diagram $G$, since
$G-f_R(G)K_0$ belongs to $\CC$, hence 
\dez{
0=
f_{R'}(G-f_{R}(G)K_0)=
f_{R'}(G)-f_R(G)f_{R'}(K_0)=f_{R'}(G)-f_R(G).
\bxx
}

\sectz{Uniqueness}

Theorem \ref{7ja09a} implies that (for each fixed $n$)
the ring $p_n(\oC\GG)\cong\oC\GG/J_{n+1}$ is finitely generated,
and that the set $\{f_R\mid R\in\MM_n\}$ of all
$n$-state partition functions
$f:\GG\to\oC$ form the affine variety $\MM_n/O_n$, since
$f_R=f_{R'}$ if and only if
$p(R)=p(R')$ for each $p\in\OO(\MM_n)^{O_n}$.
So
\dez{
\OO(\MM_n/O_n)=\OO(\MM_n)^{O_n}\cong \oC\GG/J_{n+1},
}
and the elements of $\MM_n/O_n$ are in one-to-one correspondence
with the $n$-state partition functions.

Let $\pi:\MM_n\to\MM_n/O_n$
be the corresponding projection function.
Then by Theorem \ref{23no07b}, any fiber of $\pi$ intersects $\WW_n$.
Each fiber of $\pi$ contains a unique (Zariski-)closed $O_n$-orbit
(cf.\ [8], [1]).
So the Zariski-closed $O_n$-orbits of $\MM_n$ are in one-to-one
correspondence with the $n$-state partition functions.
models.

We will show that in fact this orbit is equal to the set of elements
in the fiber that belong to $\WW_n$.
Equivalently, for each $n$-state partition function $f$ the element
$R$ of $\WW_n$ with
$f(G)=f_R(G)$ for each virtual link diagram $G$ is unique up to
the action of the orthogonal group on $R$.
So $\WW_n$ is equal to the union of the closed $O_n$-orbits.

\thm{22no08a}{
Each $R\in\WW_n$ belongs to the unique closed $O_n$-orbit $B$
contained in $\{R'\mid f_{R'}(G)=f_R(G)$ for each virtual
link diagram $G\}$.
}

\pf
Suppose not.
Then there is a polynomial $q\in\OO(\MM_n)$ with $q(B)=0$ and
$q(R)\neq 0$.
Let $U$ be the $O_n$-module spanned by $O_n\cdot q$.
The morphism $\phi:\MM_n\to U^*$ with $\phi(R')(u)=u(R')$ (for
$R'\in\MM_n$ and $u\in U$) is $O_n$-equivariant.

Let $\iota:U^*\to T(V)$ be an embedding of $U^*$ as $O_n$-submodule
of $T(V)$.
So $\iota\circ\phi$ belongs to $(\OO(\MM_n)\otimes T(V))^{O_n}=p_n(\oC\TT)$,
say $\iota\circ\phi=p_n(\tau)$ with $\tau\in\oC\TT$.
As $\phi(R)\neq 0$ (since $\phi(R)(q)=q(R)\neq 0$),
we have $p_n(\tau)(R)\neq 0$.
As $R\in\WW_n$, there is a tangle $T$ with $p_n(\tau\cdot T)(R)\neq 0$.
So $f_R(\tau\cdot T)\neq 0$.
However, for any $R'\in B$, $p_n(\tau\cdot T)(R')=
p_n(\tau)(R')\cdot p_n(T)(R')=0$,
since $p_n(\tau)(R')=\iota\circ\phi(R')=0$, as $q(B)=0$.
So $f_{R'}(\tau\cdot T)=0$ while
$f_{R}(\tau\cdot T)\neq 0$, contradicting the fact that
$f_{R'}(G)=f_{R}(G)$ for each virtual link diagram $G$.
\bx

\cor{10de07a}{
Let $R,R'\in\WW_n$ be such that $f_{R}(G)=f_{R'}(G)$ for
each virtual link diagram $G$.
Then $R'=R^U$ for some $U\in O_n$.
}

\pf
By Theorem \ref{22no08a}, $R$ and $R'$ belong to the
same $O_n$-orbit.
\bx

\sectz{Reidemeister moves}

Let $\VV_n$ be the variety of solutions of \rf{22ok07a}:
\dyz{
$\VV_n:=\{R\mid R\in\MM_n$, $R$ satisfies \rf{22ok07a}$\}$.
}
This can be translated into the following quantum tangles, corresponding
to the three Reidemeister moves:
\dyyz{
M_1:=
\mbox{\epsfxsize=0.05\hsize$\epsfcenter{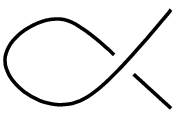}$}
-
\mbox{\epsfxsize=0.017\hsize$\epsfcenter{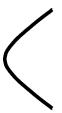}$}~,~~~
M_2:=
\mbox{\epsfxsize=0.05\hsize$\epsfcenter{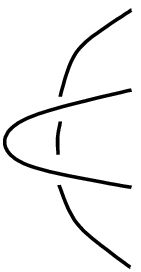}$}
-
\mbox{\epsfxsize=0.03\hsize$\epsfcenter{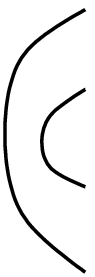}$}~,~~~
M_3:=
\mbox{\epsfxsize=0.04\hsize$\epsfcenter{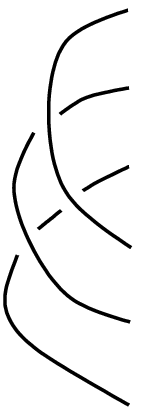}$}
-
\mbox{\epsfxsize=0.045\hsize$\epsfcenter{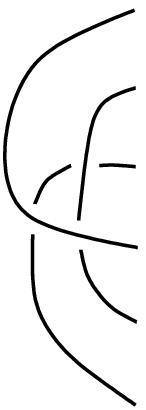}$}~.
}
Then (cf.\ \rf{22ok07a}), up to permuting tensor factors,
$f_R(M_1)=C(R)-I_n$,
$f_R(M_2)=RD(R)-I_n\otimes I_n$, and
$f_R(M_3)=E_{1,2}(R)E_{1,3}(R)E_{2,3}(R)-
E_{2,3}(R)E_{1,3}(R)E_{1,2}(R)$.

\thm{17ok07d}{
Let $f:\GG\to\oC$ and $n\in\oN$.
Then $f=f_R$ for some $R\in \VV_n$ if and only if
$f$ is an $n$-state partition function and is
invariant under Reidemeister moves.
}

\pf
Necessity is direct, by definition of $\VV_n$.
As for sufficiency,
by Theorem \ref{23no07b} we can assume that $R\in\WW_n$.
Since $f_R$ is invariant under Reidemeister moves, we have
$f_R(M_i)=0$ for $i=1,2,3$.
Hence $R\in\VV_n$.
\bx

Let again $\pi:\MM_n\to\MM_n/O_n$ be the projection function.
Then by Theorem \ref{17ok07d}, any $n$-state partition function
$f$ is a virtual link invariant if and only if the fiber $\pi^{-1}(f)$
intersects $\VV_n$.

Now $\VV_n$ is an $O_n$-invariant variety.
Then the associated affine variety $\VV_n/O_n$ is in one-to-one
correspondence with the set of $n$-state partition functions
that are invariant under Reidemeister moves.

Each fiber $\pi^{-1}(f)$ of $\pi$ contains a unique (Zariski-)closed
$O_n$-orbit.
Then $f$ is a virtual link invariant if and only
if this closed orbit is a subset of $\VV_n$.
So there is a one-to-one correspondence between the closed
$O_n$-orbits in $\VV_n$ and the Reidemeister move invariant $n$-state 
partition functions.
Let $\VV'_n$ be the set of elements $R\in \VV_n$ for which
the orbit $O_n\cdot R$ is closed.
So $\VV'_n=\VV_n\cap\WW_n$.
Then for any $R,R'\in \VV'_n$ one has: $f_R(G)=f_{R'}(G)$ for each
virtual link diagram if and only if
$R'=R^U$ for some $U\in O_n$.

\sectz{The real case}

We next give a characterization of virtual link diagram invariants
$f$ for which there exists a real $R\in\MM_n$ such that $f=f_R$.
As we will see, any real $R$ automatically belongs to $\WW_n$.

\thm{3no07b}{
Let $f$ be a real-valued virtual link diagram invariant
and let $n\in\oN$.
Then there is a real $R\in\MM_n$ with $f=f_R$
if and only if $f$ is multiplicative, $f(\det_{n+1}\cdot T)=0$
for each $2(n+1)$-tangle $T$, and the matrix
\de{10ja08h}{
(f(T\cdot T'))_{T,T'\text{ $4$\rm -tangles}}
}
is positive semidefinite.
}

\pf
Necessity follows from the fact that
\dyyz{
f_R(T\cdot T')=f_R(T)\cdot f_R(T').
}

We prove sufficiency.
As $f$ is an $n$-state partition function,
there is an algebra homomorphism $\hat f:\OO(\MM_n)^{O_n}\to\oC$
such that $\hat f\circ p_n=f$.
We must show the existence of a real-valued $R\in\MM_n$ such that
$f_{R}(G)=f(G)$ for each $G$; i.e., $p_n(G)(R)=\hat f\circ p_n(G)$ for
each $G$; in other words, by Theorem \ref{7ja09a}: $q(R)=\hat f(q)$
for each $q\in\OO(\MM_n)^{|O_n}$,

By the Procesi-Schwarz theorem [9], such an $R$ exists
if and only if
\de{9de10c}{
\hat f(dq\cdot dq)\geq 0
}
for each $q\in\OO(\MM_n)^{O_n}$ which is real-valued
on real-valued elements of $\MM_n$.
Here $d$ is the derivative function $\OO(\MM_n)\to\OO(\MM_n)\otimes\MM_n^*$.
The bilinear form $\cdot:(\OO(\MM_n)\otimes\MM_n^*)^2\to\OO(\MM_n)$
is defined by $(a\otimes b)\cdot(c\otimes d):=(b\cdot d)ac$ for
$a,c\in\OO(\MM_n)$ and $b,d\in\MM_n^*$, where $b\cdot d$ is the
standard inner product on $\MM_n^*$ (induced by the standard
inner product on $\oC^n$).

So it suffices to prove \rf{9de10c}.
We can write
\dez{
q=\sum_G\lambda_Gp_n(G),
}
with $\lambda_G\in\oR$, and only finitely many of them nonzero.
Here $G$ ranges over all virtual link diagrams.
Indeed, we can write $q=\sum_G\lambda_Gp_n(G)$ with $\lambda_G\in\oC$.
Then $q=\sum_G\kfrac12(\lambda_G+\overline\lambda_G)p_n(G)$ on the real
elements of $\MM_n$, hence on all elements of $\MM_n$.

For any virtual link diagram $G$, we define a quantum $4$-tangle $dG$.
We first define for each $v\in V(G)$, a quantum $4$-tangle $G_v$.
Let $e_1,e_2,e_3,e_4$ be the edges incident with $v$, in cyclic order,
where $e_1,e_3$ is the over-going pair.
Let $H$ be the graph obtained by deleting vertex $v$ and connecting
$e_1,\ldots,e_4$ to new vertices $v_1,\ldots,v_4$.
Then $G_v$ is defined to be half of the sum of the $4$-tangle obtained
from $H$ by giving $v_1,\ldots,v_4$ labels $1,\ldots,4$ respectively,
and the $4$-tangle obtained from $H$ by giving $v_1,\ldots,v_4$ labels
$3,4,1,2$ respectively.
Then
\dez{
dG:=\sum_{v\in V(G)}G_v.
}

Then for any $G$, $dp_n(G)=p_n(dG)$, hence for any $G,H$:
\dyyz{
dp_n(G)\cdot dp_n(H)
=
p_n(dG)\cdot p_n(dH)
=
p_n(dG\cdot dH).
}
So
\dyyz{
dq\cdot dq
=
\sum_{G,H}\lambda_G\lambda_Hdp_n(G)\cdot dp_n(H)
=
\sum_{G,H}\lambda_G\lambda_Hp_n(dG\cdot dH).
}
Therefore,
\dyyz{
\hat f(dq\cdot dq)
=
\sum_{G,H}\lambda_G\lambda_H f(dG\cdot dH)\geq 0,
}
by the positive semidefiniteness of \rf{10ja08h}.
This proves '\rf{9de10c}.
\bx

The positive semidefiniteness condition may be considered as a form
of `reflection positivity'.
The proof of the theorem is inspired by the theorem of
Szegedy [12] on edge model graph invariants.
A similar theorem (but different proof) of
Freedman, Lov\'asz, and Schrijver [2] may yield a spin model
analogue of Theorem \ref{3no07b}.

For real $R$, \rf{22ok07a} is automatically satisfied:
\dy{13de07a}{
If $R\in\MM_n$ is real, then $R\in\WW_n$.
}
Indeed, suppose $f_R(\tau)\neq 0$.
Then $f_R(\tau)\cdot f_R(\tau)\neq 0$.
Hence $f_R(\tau\cdot\tau)\neq 0$.

This implies:
\dy{17ok07b}{
Let $R\in\MM_n$ be real.
Then $f_R$ is a virtual link invariant
if and only if $R$ belongs to $\VV_n$.
}

Here it suffices to show necessity.
If $f_R$ is a virtual link invariant, then, as $R$ belongs
to $\WW_n$ by \rf{13de07a} and as $f_R$ is invariant under the
Reidemeister moves, $R$ satisfies \rf{22ok07a}.
So $R\in\VV_n$.

It is an open problem whether for any two nonisotopic virtual
link diagrams $G$ and $H$ one has $f_R(G)\neq f_R(H)$
for some $n\in\oN$ and some $R\in\MM_n$.

\section*{References}\label{REF}
{\small
\begin{itemize}{}{
\setlength{\labelwidth}{4mm}
\setlength{\parsep}{0mm}
\setlength{\itemsep}{1mm}
\setlength{\leftmargin}{5mm}
\setlength{\labelsep}{1mm}
}
\item[\mbox{\rm[1]}] M. Brion, 
Introduction to actions of algebraic groups,
{\em Les cours du C.I.R.M.} 1 (2010) 1-22.

\item[\mbox{\rm[2]}] M.H. Freedman, L. Lov\'asz, A. Schrijver, 
Reflection positivity, rank connectivity, and homomorphisms of graphs,
{\em Journal of the American Mathematical Society} 20 (2007) 37--51.

\item[\mbox{\rm[3]}] R. Goodman, N.R. Wallach, 
{\em Symmetry, Representations, and Invariants},
Springer, Dordrecht, 2009.

\item[\mbox{\rm[4]}] P. de la Harpe, V.F.R. Jones, 
Graph invariants related to statistical mechanical models:
examples and problems,
{\em Journal of Combinatorial Theory, Series B} 57 (1993) 207--227.

\item[\mbox{\rm[5]}] L. Kauffman, 
{\em Knots and Physics},
World Scientific, Singapore, 1991.

\item[\mbox{\rm[6]}] L.H. Kauffman, 
Virtual knot theory,
{\em European Journal of Combinatorics} 20 (1999) 663--690.

\item[\mbox{\rm[7]}] L.H. Kauffman, 
Introduction to virtual knot theory,
ArXiv http://arxiv.org/abs/1101.0665

\item[\mbox{\rm[8]}] H. Kraft, 
{\em Geometrische Methoden in der Invariantentheorie},
Vieweg, Braunschweig, 1984.

\item[\mbox{\rm[9]}] C. Procesi, G. Schwarz, 
Inequalities defining orbit spaces,
{\em Inventiones Mathematicae} 81 (1985) 539--554.

\item[\mbox{\rm[10]}] A. Schrijver, 
Tensor subalgebras and First Fundamental Theorems in invariant theory,
{\em Journal of Algebra} 319 (2008) 1305--1319.

\item[\mbox{\rm[11]}] A. Schrijver, 
Graph invariants in the spin model,
{\em Journal of Combinatorial Theory, Series B} 99 (2009) 502--511.  

\item[\mbox{\rm[12]}] B. Szegedy, 
Edge coloring models and reflection positivity,
{\em Journal of the American Mathematical Society}
20 (2007) 969--988.

\item[\mbox{\rm[13]}] V.G. Turaev, 
The Yang-Baxter equation and invariants of links,
{\em Inventiones Mathematicae} 92 (1988) 527--553.

\item[\mbox{\rm[14]}] H. Weyl, 
{\em The Classical Groups --- Their Invariants and Representations},
Princeton University Press, Princeton, New Jersey, 1946.

\end{itemize}
}

\end{document}